\documentclass[preprint,12pt,3p]{elsarticle}

\usepackage{float}
\usepackage{graphicx}
\usepackage{epstopdf}
\usepackage[T1]{fontenc}
\usepackage{amsmath}
\usepackage{amssymb}
\usepackage[none]{hyphenat}
\usepackage{setspace}
\usepackage{lipsum}
\makeatletter
\def\ps@pprintTitle{%
 \let\@oddhead\@empty
 \let\@evenhead\@empty
 \def\@oddfoot{}%
 \let\@evenfoot\@oddfoot}
\makeatother

\begin{document}

\begin{frontmatter}

\title{A Step-by-Step Procedure for Local Analysis of Differential Equations}

\author[add1]{Alexander Maslov}
\ead{alexander.maslov@carleton.ca}
\author[add2]{David Amundsen}
\ead{dave@math.carleton.ca}
\address[add1]{Carleton University, Department of Economics}
\address[add2]{Carleton University, School of Mathematics and Statistics}

\begin{abstract}

This note provides a detailed algorithm to the application of local (perturbation) analysis of differential equations which is normally taught at graduate math courses. Exercise books often present more abstract and simplified versions of equations for the application of perturbation techniques. The equation we study comes from the theory of competing mechanisms and describes the behavior of rational buyers in a certain environment. While in the latter literature similar equations were solved numerically, an analytical solution adds both theoretical and practical value for students and researchers.

\end{abstract}

\begin{keyword}
Ordinary Differential Equations, Perturbation Theory, Binomial Expansion
\end{keyword}

\end{frontmatter}

\doublespacing

\section{Sample Equation}
\label{sec2}

The equation used for the purpose of this paper was derived as a result of the maximization of a buyer's expected payoff in a certain environment containing an auction and a posted price. A buyer's strategy or payoff function $z:V\rightarrow R$ is mapped from his / her valuation $v$. The following economically relevant restrictions have been imposed during the construction of the model: $n\geq3$ - the number of bidders (players), $0<z(v)<p<v<1$ - where $p$ is the posted price available outside the auction. All variables and parameters belong to the set of real numbers.
The equation has the following form:
\begin{equation}\label{eq_base}
\frac{1}{2}(p-z(v))^2z(v)^{n-3}z'(v)=\frac{v^{n-2}-z(v)^{n-2}}{n-2}p-\frac{v^{n-1}-z(v)^{n-1}}{n-1}
\end{equation}

The first step is to select a point around which to perturb the equation. In this particular case such point is a point of singularity $(p,p)$, which is derived from the model\footnote{Note that it does not have to be a point of singularity for the perturbation purposes.}. Adjust the equation to the center of coordinates by subtracting $p$ from both variables, i.e. define $X=v-p$ and $Y=z-p$. Dropping the argument to simplify the notation we receive:
\begin{equation}
\frac{1}{2}(p+Y)^{n-3}Y^2Y'=\frac{(p+X)^{n-2}-(p+Y)^{n-2}}{n-2}p-\frac{(p+X)^{n-1}-(p+Y)^{n-1}}{n-1}
\end{equation}

Now, the function is not defined around $(0,0)$. Moreover $\lim_{X\to 0^{+}} Y=0$ with $X\gg X^2\gg X^3$.
The next step is to rescale the variables by introducing a small parameter $\epsilon$, which will help to keep track of their order. $\epsilon^0$ indicates leading order, $\epsilon^1$ - first perturbed term, $\epsilon^2$ - second perturbed term etc...  Setting $X=\epsilon x$ and $Y=\epsilon y$. We get:
\begin{equation} \label{eq:three}
\frac{1}{2}(p+\epsilon y)^{n-3}(\epsilon y)^2(y'+\epsilon y')=\frac{(p+\epsilon x)^{n-2}-(p+\epsilon y)^{n-2}}{n-2}p-\frac{(p+\epsilon x)^{n-1}-(p+\epsilon y)^{n-1}}{n-1}
\end{equation}

Next, apply binomial expansion to expand the polynomials according to the standard formula: 
\begin{equation}
(a+b)^n=\frac{a^n}{0!}+\frac{na^{n-1}b}{1!}+\frac{n(n-1)a^{n-2}b^2}{2!}+\frac{n(n-1)(n-2)a^{n-3}b^3}{3!}+...
\end{equation}

To derive an equation with one perturbed term we need to drop anything with $\epsilon^4$ and higher since the left part of the equation \eqref{eq:three} already contains a term $(\epsilon y)^2$, which will cancel out with corresponding powers of $\epsilon$ in the right-hand side, leaving only $\epsilon^0$ and $\epsilon^1$.
Applying the expansion to each of the right-hand side brackets (also multiplying by $p$ for the first fraction) we get:

\begin{equation}
\begin{gathered}
(p+\epsilon x)^{n-2}p=p^{n-1}+\epsilon(n-2)p^{n-2}x+\frac{1}{2}\epsilon^2(n-3)(n-2)p^{n-3}x^2+ \\ +\frac{1}{6}\epsilon^3(n-4)(n-3)(n-2)p^{n-4}x^3
\end{gathered}
\end{equation}

\begin{equation}
\begin{gathered}
(p+\epsilon y)^{n-2}p=p^{n-1}+\epsilon(n-2)p^{n-2}y+\frac{1}{2}\epsilon^2(n-3)(n-2)p^{n-3}y^2+ \\ \frac{1}{6}\epsilon^3(n-4)(n-3)(n-2)p^{n-4}y^3
\end{gathered}
\end{equation}

\begin{equation}
\begin{gathered}
(p+\epsilon x)^{n-1}=p^{n-1}+\epsilon(n-1)p^{n-2}x+\frac{1}{2}\epsilon^2(n-2)(n-1)p^{n-3}x^2+ \\ 
+\frac{1}{6}\epsilon^3(n-3)(n-2)(n-1)p^{n-4}x^3
\end{gathered}
\end{equation}

\begin{equation}
\begin{gathered}
(p+\epsilon y)^{n-1}=p^{n-1}+\epsilon(n-1)p^{n-2}y+\frac{1}{2}\epsilon^2(n-2)(n-1)p^{n-3}y^2+ \\
+\frac{1}{6}\epsilon^3(n-3)(n-2)(n-1)p^{n-4}y^3
\end{gathered}
\end{equation}

Substituting the polynomials into the right-hand side of the equation \eqref{eq:three} many terms cancel out and we get the following expression:

\begin{equation} \label{eq:nine}
\begin{gathered}
\frac{1}{n-2}\bigg(p^{n-1}+\epsilon(n-2)p^{n-2}x+\frac{1}{2}\epsilon^2(n-3)(n-2)p^{n-3}x^2+ \\
+\frac{1}{6}\epsilon^3(n-4)(n-3)(n-2)p^{n-4}x^3-p^{n-1}-\epsilon(n-2)p^{n-2}y- \\
\frac{1}{2}\epsilon^2(n-3)(n-2)p^{n-3}y^2-\frac{1}{6}\epsilon^3(n-4)(n-3)(n-2)p^{n-4}y^3\bigg)- \\
-\frac{1}{n-1}\bigg(p^{n-1}+\epsilon(n-1)p^{n-2}x+\frac{1}{2}\epsilon^2(n-2)(n-1)p^{n-3}x^2+ \\
\frac{1}{6}\epsilon^3(n-3)(n-2)(n-1)p^{n-4}x^3-p^{n-1}-\epsilon(n-1)p^{n-2}y- \\
-\frac{1}{2}\epsilon^2(n-2)(n-1)p^{n-3}y^2-\frac{1}{6}\epsilon^3(n-3)(n-2)(n-1)p^{n-4}y^3\bigg)=\\
=\frac{1}{2}\bigg(\epsilon^2p^{n-3}(y^2-x^2)+\epsilon^3\frac{2(n-3)}{3}p^{n-4}(y^3-x^3)\bigg)
\end{gathered}
\end{equation}

Now, substitute the final expression of the equation \eqref{eq:nine} into the right-hand side of the equation \eqref{eq:three}, dividing both parts by $\frac{1}{2}(\epsilon y)^2(p+\epsilon y)^{n-3}$ and expanding $(p+\epsilon y)^{n-3}$ in the denominator of the right-hand side:

\begin{equation} \label{eq:ten}
y'+\epsilon y'=\frac{\epsilon^2p^{n-3}(y^2-x^2)+\epsilon^3\frac{2(n-3)}{3}p^{n-4}(y^3-x^3)}{(p^{n-3}+\epsilon(n-3)p^{n-4}y) (\epsilon y)^2}
\end{equation}

It is readily seen that $\epsilon^2$ cancel out. Since $n$ is a fixed number of players, around the singularity point it is true that $n\gg y$. Therefore, we can apply geometric series transformation $\frac{1}{1-q}=1+q+q^2+q^3...$ to the polynomial in the denominator (note that the order of $\epsilon$ has to be consistent, and since the latter is multiplied by $\epsilon^2$, the expanded order coincides with the one truncated in the numerator):
\begin{equation}
\frac{1}{p^{n-3}}\frac{1}{(1-\epsilon\frac{3-n}{p}y)}=\frac{1}{p^{n-3}}(1-\frac{\epsilon(n-3)y}{p})
\end{equation}

Note, that we truncate the series after the second term, because the third term already includes $\epsilon^2$, which is higher than the order of $\epsilon$ received for equation \eqref{eq:ten}. Therefore, equation \eqref{eq:ten} transforms into:
\begin{equation} \label{eq:twelve}
y'+\epsilon y'=\frac{p^{n-3}(y^2-x^2)+\epsilon\frac{2(n-3)}{3}p^{n-4}(y^3-x^3)}{p^{n-3}y^2} \bigg(1-\frac{\epsilon(n-3)y}{p}\bigg)
\end{equation}

Since we have already canceled out $\epsilon^2$ in both parts of the equation, the newly truncated order is just $\epsilon$, so we drop any term with the power of $\epsilon$ greater than $1$. Hence, the right-hand side of the equation \eqref{eq:twelve} becomes:

\begin{equation}
\begin{gathered}
\frac{y^2-x^2}{y^2}+\epsilon \frac{2(n-3)(y^3-x^3)}{3py^2} \bigg(1-\frac{\epsilon(n-3)y}{p}\bigg)= \\
\frac{y^2-x^2}{y^2}+\epsilon \frac{2(n-3)(y^3-x^3)}{3py^2}-\epsilon \frac{y^2-x^2}{y^2}\frac{(n-3)y}{p}-\epsilon^2 \frac{2(n-3)^2y(y^3-x^3)}{3p^2y}= \\
\frac{y^2-x^2}{y^2}+\epsilon \frac{(3-n)(y-x)^2(2x+y)}{3py^2}
\end{gathered}
\end{equation}

Note, that the term with $\epsilon^2$ has been dropped and each order of the variable(s) corresponds to the order of $\epsilon$. Thus, we have transformed equation \eqref{eq:three} into an equation consisting of a leading order and one perturbed term\footnote{To receive more perturbed terms add more orders of binomial expansion in the beginning.}:

\begin{equation} \label{eq:fourteenth}
\begin{gathered}
y_0'+\epsilon y_1'=\frac{y^2-x^2}{y^2}+\epsilon \frac{(3-n)(y-x)^2(2x+y)}{3py^2}
\end{gathered}
\end{equation}

\section{Finding the Leading Order}
\label{sec3}

The leading order is the term with $\epsilon^0$, i.e. (capital letters here are used solely to emphasise the leading nature of the equation):
\begin{equation} \label{eq:fifteenth}
\frac{dY}{dX}=\frac{Y^2-X^2}{Y^2}
\end{equation}

The equation above does not have a unique solution, but can be solved approximately by the method of dominant balance, geometric approximation or by a standard approach to homogeneous ODE, which we will proceed with here:

Replace $Y=XF(X)$ to get:
\begin{equation}
\begin{gathered}
X\frac{dF}{dX}+F=\frac{F^2-1}{F^2} \ \rightarrow \ X\frac{dF}{dX}=\frac{F^2-1-F^3}{F^2} \\
\frac{dX}{X}=\frac{F^2}{-F^3+F^2-1}dF \ \rightarrow \ ln(X)=\int\frac{F^2}{-F^3+F^2-1}dF= \\
\int\frac{F^2}{-(F-r_1)(F-r_2)(F-r_3)}dF=-\frac{r_1^2\ln\left|{F-r_1}\right|}{(r_1-r_2)(r_1-r_3)}- \\ -\frac{r_2^2\ln\left|{F-r_2}\right|}{(r_2-r_1)(r_2-r_3)}-\frac{r_3^2\ln\left|{F-r_3}\right|}{(r_3-r_1)(r_3-r_2)}+c \\ X=C\  exp\Bigl(-\frac{r_1^2\ln\left|{F-r_1}\right|}{(r_1-r_2)(r_1-r_3)} -\frac{r_2^2\ln\left|{F-r_2}\right|}{(r_2-r_1)(r_2-r_3)}-\frac{r_3^2\ln\left|{F-r_3}\right|}{(r_3-r_1)(r_3-r_2)}\Bigr)\\
X=C (\left|{F-r_1}\right|)^{\frac{-r_1^2}{(r_1-r_2)(r_1-r_3)}}(\left|{F-r_2}\right|)^{\frac{-r_2^2}{(r_2-r_1)(r_2-r_3)}}(\left|{F-r_3}\right|)^{\frac{-r_3^2}{(r_3-r_1)(r_3-r_2)}}
\end{gathered}
\end{equation}

Note that in the solution above $F\ne r_1$, $F\ne r_2$ and $F\ne r_3$ because we automatically assumed that these values are not zero when dividing by the polynomial. From the parametric expression it is seen that when $F$ approaches either of the roots the value of the function becomes infinitely large, which may not be a solution. Therefore, let us check whether the roots of the polynomial are the solutions to the equation \eqref{eq:fifteenth}. In other words we check if $Y=kX$ solves this equation.

$$k=\frac{(kX)^2-X^2}{(kX)^2} \Rightarrow k^3=k^2-1$$

The last expression is exactly the equation which roots we have already found. Being only interested in the real root, the solution is thus $k=-0.7549$ or $Y=-0.7549X$, which can be verified by substituting it back into the equation \eqref{eq:fifteenth}.

\section{Searching for the Perturbed Term}
\label{sec4}

Having found the leading order, we can solve for a correction or perturbed term. To do this we need to balance the terms on both sides, i.e. 
\begin{equation}
\epsilon^0y_0'+\epsilon^1y_1'+...+\epsilon^ny_n'=\epsilon^0(.)+\epsilon^1(.)+...\epsilon^n(.)
\end{equation}

In section \ref{sec3} we have already found the solution to the leading term ODE $y_0'=(.)$. Now, plug the symbolic conceived solution including the perturbed term into the equation \eqref{eq:fourteenth}:
\begin{equation} \label{eq_sixteenth}
y_0'+\epsilon y_1'=\frac{(y_0+\epsilon y_1)^2-x^2}{(y_0+\epsilon y_1)^2}+\epsilon\frac{(3-n)(2x+y_0+\epsilon y_1)(y_0+\epsilon y_1-x)^2}{3p(y_0+\epsilon y_1)^2}
\end{equation}

Next, move the denominator to the numerator by the same way as in section \ref{sec2}\footnote{Note, that here it is squared, therefore $(\frac{1}{1-q})^2=(1+q+q^2+q^3...)^2$.}. Within the set limit for $\epsilon$ we are only interested in the first two elements of the square of the sum of the first two terms: $(1+2q$), because all other terms will have higher orders. Hence, the denominator of the right-hand side can be rewritten as:
\begin{equation}
\frac{1}{(y_0+\epsilon y_1)^2}=\frac{1}{y_0^2}\bigg(1-\frac{2\epsilon y_1}{y_0}\bigg)
\end{equation}

Multiplying the transformed denominator on the equation \eqref{eq_sixteenth} and dropping terms with $\epsilon$ higher than the order one we get:
\begin{equation}
\begin{gathered}
\bigg(1-\frac{2\epsilon  y_1}{y_0}\bigg)\bigg(\frac{(-x^2+(y_0+\epsilon y_1)^2)}{y_0^2}-\frac{\epsilon(n-3)(-x+y_0+\epsilon y_1)^2(2x+y_0+\epsilon y_1)}{3y_0^2p}\bigg)= \\
=\frac{y_0^2-x^2}{y_0^2}+\epsilon \frac{1}{3y_0^2}\bigg(\frac{(3-n)y_0^3+3(n-3)x^2y_0-2(n-3)x^3+}{p}+\frac{6x^2y_1}{y_0}\bigg)+ \\
+\epsilon^2 \frac{1}{3y_0^2}\bigg(\frac{4(n-3)x^3y_1}{py_0}+\frac{(3-n)y_0^2y_1}{p}-\frac{3(n-3)x^2y_1}{p}-9y_1^2\bigg)+ \\
+\epsilon^3 \frac{1}{3y_0^2}\bigg(\frac{3(n-3)y_0y_1^2}{p}-\frac{6(n-3)x^2y_1^2}{py_0}-\frac{6y_1^3}{y_0}\bigg)+
\epsilon^4\frac{5(n-3)y_1^3}{3py_0^2}+\epsilon^5\frac{2(n-3)y_1^4}{3py_0^3}=\\
=\frac{y_0^2-x^2}{y_0^2}+\epsilon\frac{6x^3y_0-2nx^3y_0-9x^2y_0^2+3nx^2y_0^2+3y_0^4-ny_0^4+6px^2y_1}{3py_0^3}
\end{gathered}
\end{equation}

Thus, the ODE for the perturbed term is (note, that it is linear in $y_1$):
\begin{equation}
y_1'=\frac{6x^3y_0-2nx^3y_0-9x^2y_0^2+3nx^2y_0^2+3y_0^4-ny_0^4+6px^2y_1}{3py_0^3}
\end{equation}

Substituting $y_0=-0.754878x$ into the equation above we get:
\begin{equation}
y_1'=-\frac{0.775(-8.684x^4+2.895nx^4+6px^2y_1)}{px^3}
\end{equation}

This linear ODE can be easily solved for:
\begin{equation}
y_1=\frac{(1.012-0.3373n)x^2}{p}
\end{equation}

Adding both terms together we receive the final approximate solution for the equation \eqref{eq:three}:
\begin{equation}
y=-0.755x+\frac{(1.012-0.3373n)x^2}{p}
\end{equation}

Note, that this equation is adjusted to the center of coordinates, it has negative slope ($n\geq 3$) and each power of a variable maps into the corresponding order of the true polynomial. 

Substituting back $Y=z(v)-p$ and $X=v-p$ and setting $r=1.012-0.3373n$ to simplify the expression produces the final equation:

\begin{equation}\label{eq_seventeenth}
z(v)=(\frac{7}{4}+r)p-(\frac{3}{4}+2r)v+\frac{rv^2}{p}
\end{equation}

Lastly, we would like to verify how well equation \eqref{eq_seventeenth} describes the original equation. Equation \eqref{eq_base} is a $n$-order differential equation which has an infinite number of solutions. We are interested in the ones which pass through the initial condition. Recall, however, that the initial condition, i.e. point ($p,p$), is not defined, so we will use the point ($p+\epsilon,p-\epsilon$) to test it, where $\epsilon$ is robust to perturbations. The results are presented on Figure 1.

As expected, with an increase in $n$ the agreement dissipates, which is natural when a quadratic (all other terms are set to zero, i.e. $0\times x^3$, $0\times x^4$ etc..) function approximates $n$-order nonlinear differential equation. Introduction of additional perturbed terms will allow the approximate solution to describe the original equation with increased precision.

\section{Conclusion}
\label{sec5}

In this paper we have shown how methods of asymptotic analysis of ODE, in particular - perturbation theory, can be used to find closed form solutions for high-order nonlinear ODE. While differential equations may differ from one discipline to another, the algorithm of analysis remains the same, requiring a researcher to choose a point, around which the function can be expanded using any type of series. In the economic context, for example, it is usually substantiated by the assumptions and consequent logic of the constructed model. Local analysis plays an important role allowing to establish approximate functional relationship, which can further be used in relevant economic and other theories. Finally, the procedure to find such a solution is a great learning tool for graduate students studying local analysis.

\newpage
\begin{figure}
  \centering
  \caption{Simulating the degree of agreement between the approximate solution and original differential equation for different parameters of $p$ and $n$.}
   \centerline{\includegraphics[width=7 in, height=8 in]{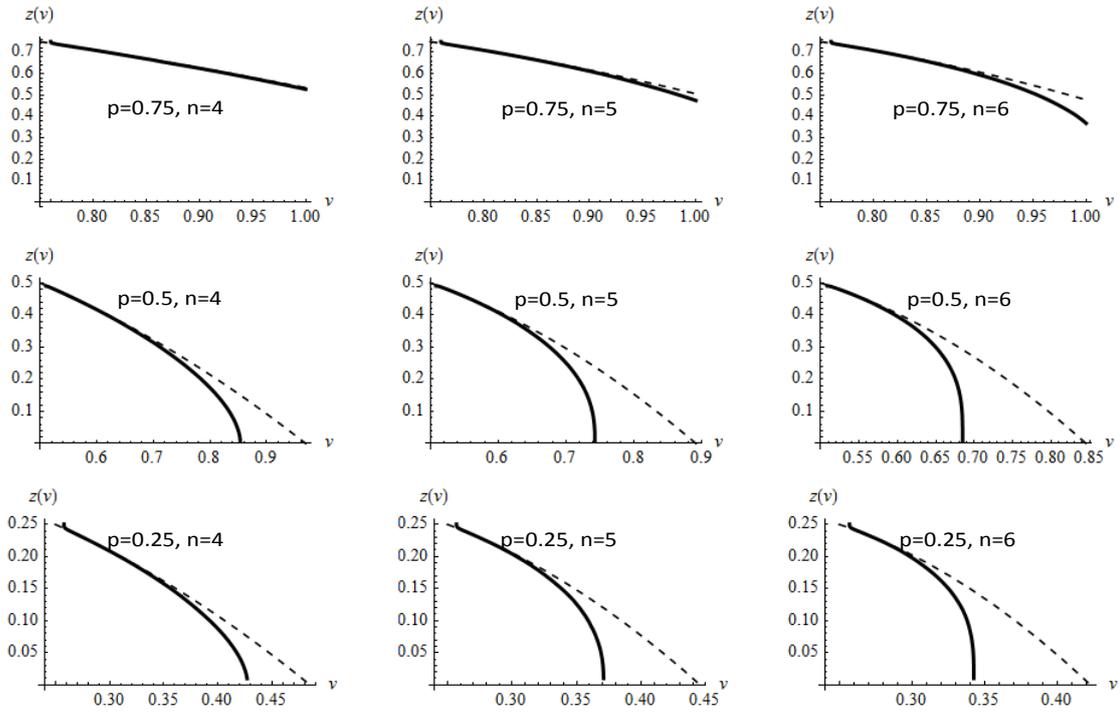}}
\end{figure}

\end{document}